%% file: root.tex
\def\mytitle{Vehicle mission guidance by symbolic optimal control}
\def\myname{Alexander Weber, Florian Fiege and Alexander Knoll}
\def\tsp{travelling salesman problem}
\def\Tsp{Travelling salesman problem}
\def\Cvrp{Capacitated vehicle routing problem}

\def\mykeywords{\Cvrp, \Tsp, symbolic optimal control, abstraction-based synthesis}
\def\arxiv{}
\documentclass[letterpaper, 10pt, conference]{ieeeconf}
\IEEEoverridecommandlockouts
\usepackage{cite}
\usepackage[english]{babel}
\usepackage{amsmath}
\usepackage{amssymb}
\usepackage{bm}

\usepackage{amsthm}
\usepackage{stmaryrd}
\usepackage{tikz}
\usetikzlibrary{calc,shapes,arrows,automata}
\usepackage[]{caption}
\usepackage{algorithm}
\usepackage{algpseudocode}
\algnewcommand\algorithmicinput{\textbf{Input:}}
\algnewcommand\algorithmicoutput{\textbf{Output:}}
\algnewcommand\algorithmicparameter{\textbf{Parameter:}}
\algnewcommand\Input{\item[\algorithmicinput]}
\algnewcommand\Output{\item[\algorithmicoutput]}
\algnewcommand\Parameter{\item[\algorithmicparameter]}
\swapnumbers
\newtheoremstyle{rem}{\topsep}{\topsep}{\normalfont}{0pt}{\bfseries}{.}{ }{\thmname{#1}\thmnumber{#2}\thmnote{ \textup{(#3)}}}
\newtheorem{definition}{Definition}[section]

\theoremstyle{rem}
\newtheorem{example}[definition]{\pushQED{\qed}Example}

\makeatletter
\def\endexample{\popQED\@endtheorem}
\def\term#1%
{\@nomath\term\ifdim\fontdimen\@ne\font>\z@%
\textbf{#1}\else\textit{#1}\fi}
\def\intcc#1{\ensuremath{\left[#1\right]}}
\def\intco#1{\ensuremath{\left[#1\right[}}
\def\intoc#1{\ensuremath{\left]#1\right]}}
\makeatother
\def\powerset#1%
{\mathcal{P}(#1)%
}
\newcommand{\R}{\mathbb{R}}

\def\implies{\relax\ifmmode\mathrel{\Rightarrow}\else$\implies$ \fi}
\DeclareMathOperator*{\argmin}{arg\,min}
\DeclareMathOperator*{\len}{len}
\usepackage{paralist}
\usepackage{hyperref}
\hypersetup{
colorlinks=true,%
breaklinks=true,%
pdfdisplaydoctitle=true,%
linkcolor={black},%
citecolor={black},
urlcolor={blue},
pdfstartview={FitV},%
pdftitle={\mytitle},
pdfauthor={\myname},%
pdfsubject={Submitted to arxiv on \today.},
pdfkeywords={\mykeywords}
}
\usepackage{subcaption}
\graphicspath{{figures/}} %
\ifx\arxiv\undefined
\else
\usepackage{fancyhdr}
\fi
\title{\bf \LARGE \mytitle}
\author{
\myname
\thanks{
The authors are with the
Munich University of Applied Sciences,
Dept. of Mechanical, Automotive and Aeronautical Eng.,
80335 M\"unchen, Germany.
}
\thanks{This work has been supported by the German Federal Ministry of Education and Research (Project ARCUS; No. 13FH734IX6). %
}%
}%
\begin{document}
\maketitle

\begin{abstract}
Symbolic optimal control is a powerful method 
to synthesize algorithmically correct-by-design 
state-feedback controllers 
for nonlinear plants. 
Its solutions are (near-)optimal 
with respect to a given cost function.
In this note, it is demonstrated 
how symbolic optimal control can be used to 
calculate controllers 
for an optimized routing 
guidance of vehicle systems
in continuous state space. 
In fact, the capacitated vehicle routing problem and 
a variant of travelling salesman problem are investigated. 
The latter problem
has a relevant application 
in case of loss of vehicles during mission.
A goods delivery scenario and 
a reconnaissance mission,
involving bicycle and aircraft dynamics respectively, 
are provided as examples.

\end{abstract}
\section{Introduction}
\label{s:introduction}
Vehicle routing problems exist in numerous variations and 
often reflect relevant problems in applied logistics 
\cite{DantzigRamser59,TothVigo02}. 
The raw vehicle routing problem (VRP) 
consists of a weighted graph and 
asks for a set of routes for a fleet of vehicles 
for visiting all ``customer" nodes optimally 
when starting and returning to the ``depot" node.
Techniques of combinatorial optimization provide solvers for it 
which return either exact or approximately optimal solutions. 
Mostly, vehicle routing problems are posed 
on ordinary directed or undirected graphs. Thus, they 
neglect any kind of dynamics behind the routes.
Some exceptions consider Dubins vehicle dynamics and 
a fixed variant of the VRP, e.g.
\cite{AndersonMilutinovic13,
ManyamRathinamDarbha15,
MansouriLagriffoulPecora17}. 
In our preceding work \cite{WeberKnoll21}, 
we considered the classical \tsp{} for arbitrary 
nonlinear dynamics.

This note follows the approach in our previous work, 
i.e. assumes nonlinear vehicle dynamics in discrete time and 
continuous state space. In fact, the dynamics
are given by
\begin{equation}
\label{e:discrete:dynamics}
x(t+1) \in F(x(t),u(t))
\end{equation} 
with state signal $x$ in $\mathbb{R}^n$, 
input signal $u$ and a strict set-valued map $F$. 
Customer locations as well as the depot are assumed 
by non-empty sets in $\mathbb{R}^n$ and the cost for the mission 
is obtained by accumulating a non-negative cost term 
along the routes.
In our example sections, 
we consider vehicle dynamics 
which are given in continuous time and 
are subject to additive disturbances 
bounded in $W \subseteq \mathbb{R}^n$. 
In fact, the dynamics is 
given by the differential inclusion
\begin{equation}
\label{e:cont:dynamics}
\dot x \in f(x,u) + W
\end{equation}
where $f$ is an ordinary map. 
By sampling, \eqref{e:cont:dynamics} is transferred to \eqref{e:discrete:dynamics} 
making our methods applicable to dynamics \eqref{e:cont:dynamics}.

The focus of this work is on two variants of the VRP, 
namely the capacitated vehicle routing problem (CVRP) 
and a dynamic variant of the \tsp{} (TSP). 
The algorithms to present will provably solve 
previous problems with respect to the qualitative task, 
i.e. the control task of visiting a prescribed set 
of targets under some further constraints on the order of visit. 
(In this work, we refer to this 
special coverage specification as ``mission".) 
In addition, we provide novel heuristics within 
the algorithms to quantitatively reduce the cost 
for the mission in total.

As we demonstrate, 
our methods can be used 
for general vehicle mission guidance 
and task reassignments. 
As an example for the latter, 
we consider reconnaissance missions 
with uninhabited aerial vehicles (UAVs) 
where at some point during mission
a UAV suffers a damage. 
In this case, 
the remaining UAVs must take over 
the tasks of the broken UAV. 
This kind of 
scenario is not new, e.g. is investigated in
\cite{BellinghamTillersonRichardsHow03,
MiaoZhongYinZouLuo17,
FuMaoHeYuXie19}. 
In these works, nonlinear vehicle dynamics are 
basically neglected and
solutions based on 
infinite digraphs are presented. 
A drawback of our approach is that 
our methodology does not ensure 
collision avoidance if the mission 
involves more than one vehicle.
To summarize, 
our novel methods solve coverage problems 
under certain additional constraints
that involve nonlinear dynamics, 
hard state or input constraints, 
uncertainties or measurement errors. 
These features come out-of-the-box with
correct-by-design controllers which 
have been heuristically optimized 
to reduce mission costs. 
To this end, 
the framework of symbolic optimal control \cite{ReissigRungger18} 
and the related computational means \cite{ReissigWeberRungger17} 
are used.

The remaining part of this note is 
organized as follows. 
After introducing notation (Section \ref{s:notation}) 
a review on symbolic optimal control is given (Section \ref{ss:symbolicoptimalcontrol})
and some known results on coverage problems are stated (Section \ref{ss:coverage problems}). 
Sections \ref{s:vrp} and \ref{s:tsp} contain the novel methods on 
the CVRP and the dynamic TSP, respectively. 
Each of the two previous sections contains 
a subsection with an experimental evaluation. 
Conclusions are given in Section \ref{s:conclusions}.
\ifx\arxiv\undefined
\else
\thispagestyle{fancy}
\fi
\section{Notation}
\label{s:notation}
The symbols 
$\mathbb{R}$, 
$\mathbb{R}_+$, 
$\mathbb{Z}_+$, 
$\mathbb{N}$ stand for 
the set of real numbers, 
the subset of non-negative real numbers, 
the set of non-negative integers and
the set of positive integers, respectively. 
The closed and half-open intervals in $\mathbb{R}$ are denoted by 
$\intcc{a,b}$, $\intco{a,b}$, $\intoc{a,b}$ for $a,b \in \mathbb{R}$. 
For discrete intervals, we write 
$\intcc{a;b}$, 
$\intco{a;b}$, 
$\intoc{a;b}$, e.g.
$\intoc{a;b} = \intoc{a,b} \cap \mathbb{Z}$.
If $A,B$ are sets 
the notation $f \colon A \rightrightarrows B$ means 
that $f$ is a set-valued map. 
$f$ is \emph{strict} if 
$f(a) \neq \emptyset$ for all $a \in A$, 
where $\emptyset$ stands for the empty set. 
If $f \colon A \times B \to C$ is an ordinary map then
$f(\cdot,b)$ stands for the map $A \to C$ mapping $a$ to $f(a,b)$. 
If $f$, $g$ are maps then $f \circ g$ stands for the composition of $f$ and $g$.
$A^B$ denotes the set of all (ordinary) maps $B \to A$. 
E.g. $A^{\mathbb{Z}_+}$ stands for 
the set of all sequences $(a_0,a_1,\ldots)$ with 
$a_i \in A$ for all $i \in \mathbb{Z}_+$. Elements of $A^{\mathbb{Z}_+}$ are called \term{signals}.
%
%
The notation $(a_i)_{i=1}^k$ stands for the finite sequence $(a_1,\ldots,a_k)$. 
The length of a finite sequence $a$
is denoted by $\len a$, e.g. $\len (a_1,a_2,a_3) = 3$.
%
%
%
%
%
\section{Review: Symbolic optimal control and quantitative coverage problems}
\label{s:review}
We review the concept of symbolic optimal control
as established in \cite{ReissigRungger18,WeberKnoll20} 
since the novel methods are built on top.
Our methods also use basic techniques for solving
quantitative reachability problems or more generally 
coverage problems. 
Therefore, a review on these specifications is also included.
\subsection{Symbolic optimal control}
\label{ss:symbolicoptimalcontrol}
\begin{figure}
\centering
\input{figures/highlevelloop.tikz}
\caption{\label{fig:controller}Controller structure \cite{WeberKnoll21} used in this work. 
The scheduler switches between memoryless controllers 
$\mu_1,\ldots,\mu_K \colon X \rightrightarrows U \times \{0,1\}$. 
Switching is triggered when the stopping signal $v$ 
of the active controller (gray) changes from $0$ to $1$.}
\vspace*{-\baselineskip}
\end{figure}
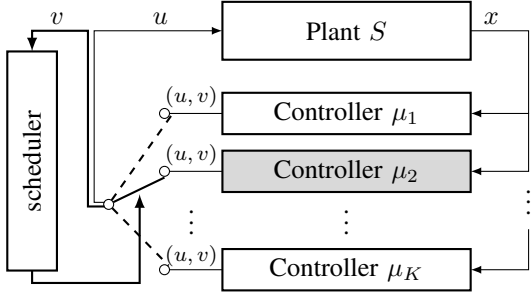

Our methods apply to dynamical systems that 
can be cast to \term{transition systems}.
By the latter we mean a triple
\begin{equation}
\label{e:system}
(X,U,F)
\end{equation}
which defines dynamics \eqref{e:discrete:dynamics}, where 
$F \colon X \times U \rightrightarrows X$ is a strict set-valued map.
The non-empty sets $X$ and $U$ are called 
\term{state} and 
\term{input space}, respectively. 
Let $S$ denote \eqref{e:system}. 
A \term{controller} for $S$ is a strict set-valued map
\begin{equation}
\label{e:controller}
\mu \colon \bigcup_{T \in \mathbb{Z}_+}\nolimits X^{\intcc{0;T}} \rightrightarrows U \times \{0,1\}
\end{equation}
and shall influence by concept the behaviour of $S$. 
The second component of the image of $\mu$
reports if the controller is indeed enabled ('$0$') or 
disabled ('$1$') \cite{ReissigRungger13}. 
In this work, controller \eqref{e:controller}
will be composed of several memoryless controllers 
$X \rightrightarrows U \times \{0,1\}$ 
as shown in Fig.~\ref{fig:controller}.

The \term{closed-loop behaviour} (initialized at $p \in X$) is the 
set of all signals $(u,v,x) \in (U \times \{0,1\} \times X)^{\mathbb{Z}_+}$
that satisfy 
\eqref{e:discrete:dynamics} and 
$(u(t),v(t)) \in \mu(x|_{\intcc{0;t}})$ 
for all $t \in \mathbb{Z}_+$. 
This set is denoted by $\mathcal{B}_p(\mu \times S)$.
In this work, the operation of the closed loop
causes costs and the accomplishment of a control task
in \emph{finite} time is mandatory. In fact, 
the cost for operating the closed loop 
is given by virtue of the cost functional 
\begin{subequations}
\label{e:totalcost}
\begin{equation}
\label{e:costfunctional:declaration}
J \colon (U \times \{0,1\} \times X)^{\mathbb{Z}_+} \to \mathbb{R}_+ \cup \{\infty\}.
\end{equation}
It is defined by $J(u,v,x) = \infty$ for perpetual operation, 
i.e. if the \term{stopping signal} $v$ is identically zero, 
and otherwise by
\begin{equation}
\label{e:costfunctional:definition}
J(u,v,x) = G(x|_{\intcc{0;T}}) + \sum_{t = 0}^{T-1} g(x(t),x(t+1),u(t))
\end{equation}
\end{subequations}
where $T = \inf v^{-1}(1) < \infty$. 
In \eqref{e:costfunctional:definition}, the \term{trajectory cost}
\begin{equation}
\label{e:trajectorycost}
G \colon \bigcup_{T \in \mathbb{Z}_+}\nolimits X^{\intcc{0;T}} \to \mathbb{R}_+ \cup \{ \infty \}
\end{equation}
rates the trajectory until stopping, while the \term{running cost}
\begin{equation}
\label{e:runningcost}
g \colon X \times X \times U \to \mathbb{R}_+ \cup \{ \infty \}
\end{equation}
is accumulative and takes also into account the input signal. 
All in all, an optimal control problem 
is as follows\cite{WeberKnoll20,ReissigRungger18}.
\begin{definition}
Let $(X,U,F)$ be a system and 
let $G$ and $g$ be as in \eqref{e:trajectorycost} and \eqref{e:runningcost}. 
An \term{optimal control problem} is a 5-tuple
\begin{equation}
\label{e:ocp}
(X, U, F, G, g).
\end{equation}
\end{definition}
The goal of symbolic optimal control is 
to find a controller so that the operation
of the closed loop is at \emph{finite} cost and ideally 
minimized for the worst case trajectory. 
To formalize this goal, 
we review the notion of \emph{closed-loop performance}
and \emph{value function}. 
Subsequently, $\mathcal{F}(X,U)$ denotes the 
set of all controllers of the form \eqref{e:controller}.
\begin{definition}
\label{def:performancefunction}
Let $\Pi$ be the optimal control problem in \eqref{e:ocp}. 
Let $J$ be the cost functional \eqref{e:totalcost} as defined for $\Pi$. 
The map $L \colon X \times \mathcal{F}(X,U) \to \mathbb{R}_+ \cup \{ \infty \}$
defined by 
\begin{equation*}
L(p,\mu) := \sup_{(u,v,x) \in \mathcal{B}_p(\mu \times S) } J(u,v,x)
\end{equation*}
is called \term{performance function} of $\Pi$. 
The \term{closed-loop performance} of $\mu \in \mathcal{F}(X,U)$ is the 
function $L(\cdot, \mu)$.
\end{definition}
As usually in optimal control, 
the value function is related to 
the best possible closed-loop performance.
Our methods will actually return suboptimal solutions, 
which are also defined next 
and turn out to be useful in applications.
\begin{definition}
\label{def:valuefunction}
Let $\Pi$ and $L$ be as in Definition \ref{def:performancefunction}.
The \term{value function} of $\Pi$ is the map assigning $p \in X$ to
\begin{equation}
V(p) := \inf_{\nu \in \mathcal{F}(X,U)} L(p,\nu).
\end{equation}
Let $A \subseteq X$ and let $\mu \in \mathcal{F}(X,U)$.
The controller $\mu$ 
\term{solves $\Pi$ suboptimally on $A$} 
if $V(p) < \infty$ implies $L(p,\mu) < \infty$ for all $p \in A$.
$\mu$ is \term{optimal} if $V = L(\cdot,\mu)$ and 
the pair $(V,\mu)$ is an \term{optimal solution} of $\Pi$. 
\end{definition}
\subsection{Quantitative reach-avoid and coverage problems}
\label{ss:coverage problems}
In the special case that $\Pi$ in Def.~\ref{def:valuefunction} 
is a quantitative reach-avoid problem and $X$ and $U$ are finite,
computational techniques are available to solve $\Pi$. 
These techniques
are the algorithmic core of the novel methods to present
and briefly described subsequently.

A quantitative reach-avoid problem
asks to steer the state of the system 
to a specified target set in finite time
while minimizing the cost functional $J$ in \eqref{e:totalcost}. 
The constraint of obstacle avoidance can be 
formulated using the running cost \eqref{e:runningcost} 
\cite[Ex.~III.5]{ReissigRungger18}.
The precise definition is:
\begin{definition}[\!\!\cite{WeberKnoll20}]
\label{def:reachavoid}
Let $\Pi$ be an optimal control problem 
of the form \eqref{e:ocp} such that $G$ is
defined by 
\begin{equation*}
G(x|_{\intcc{0;t}}) = \begin{cases}
G_0(x(t)), & \text{if } x(t) \in A \\
\infty, & \text{otherwise}
\end{cases}
\end{equation*}
where $A \subseteq X$ non-empty and 
$G_0 \colon X \to \mathbb{R}_+ \cup \{\infty\}$. 
Then $\Pi$ is called \term{(quantitative) reach-avoid problem} 
associated with $S$, $g$, $A$ and $G_0$, where
$S$ denotes the system $(X,U,F)$.
\end{definition}
The parameter function $G_0$ above is significant \cite{WeberKnoll20}
in order to achieve partial optimality for 
the type of coverage problems considered in this work. 
The next definition generalizes Def.~\ref{def:reachavoid}. 
In fact, coverage asks to steer the system state 
to a given list of target sets in finite time (in any order). 
\begin{definition}
\label{def:coverage}
Let $S$ and $\Pi$ be as in Def.~\ref{def:reachavoid}
such that $G$ is defined by
\begin{equation*}
G(x|_{\intcc{0;t}}) = \begin{cases} G_0(x|_{\intcc{0;t}}), & \text{if } \ \forall_{i \in \intcc{1;N}} \exists_{s \in \intcc{0;t}} : x(s) \in A_i\\
\infty, & \text{otherwise} 
\end{cases}
\end{equation*}
where $A_1,\ldots,A_N \subseteq X$ are non-empty and $G_0 \colon \cup_{T \in \mathbb{Z}_+} X^{\intcc{0;T}} \to \mathbb{R}_+ \cup \{\infty\}$.
Then $\Pi$ is called \term{(quantitative) coverage problem} associated with 
$S$, $g$, $A_1,\ldots,A_N$ and $G_0$.\looseness=-1
\end{definition}
By means of $G_0$ further constraints 
to a coverage problem on the trajectory can be added. 
A coverage problem associated 
with $G_0$ defined as identically $0$ 
can be solved by 
the fixed-point iteration algorithm 
shown in Fig.~\ref{a:coverage}. 
\begin{figure}[b]
\begin{algorithmic}[1]
\Function{SolveCoverage}{$S,g,A_1,\ldots,A_N$}
\State{\label{alg:queue}$Q \gets \{ 1,\ldots, N \}$\Comment{\small{}Initialize a queue}}
\State{$(A_1',\ldots,A_N') \gets (A_1,\ldots,A_N)$\Comment{\small{}Subsets of targets}}
\While{\label{alg:while}$Q\neq \emptyset$}
\State{Pick $i \in Q$}
\State{$Q \gets Q \setminus \{i\}$}
\State{\label{alg:valuefunction}$(V_i,\mu_i) \gets \textsc{SolveReach}(S,g,A_i',0)$}
\ForAll{$j \in \{1,\ldots,N\} \setminus \{i\}$}
\State{$A_\mathrm{tmp} \gets A'_j \setminus V_i^{-1}(\infty)$\Comment{\small{}$A_\mathrm{tmp} \subseteq A'_j$}}
\If{$A_\mathrm{tmp} = \emptyset$}\Comment{\small{}Solving failed}
\State{\Return{``Problem can't be solved"}}
\ElsIf{$A'_j \neq A_\mathrm{tmp}$}
\State{$A'_j \gets A_\mathrm{tmp}$\Comment{\small{}Shrink $A'_j$}}
\State{$Q \gets Q \cup \{j\}$\Comment{\small{}Solve again for $j$}}
\EndIf{}
\EndFor{}
\EndWhile{\label{alg:while:end}}
\State{\Return $((V_1,A_1',\mu_1),\ldots,(V_N,A_N',\mu_N))$}
\EndFunction
\end{algorithmic}
\caption{\label{a:coverage}Algorithm to solve 
the coverage problem associated with $S$, $g$, $A_1,\ldots,A_N$ and the zero function \cite{WeberKnoll21}.}
\vspace*{-\baselineskip}
\end{figure}
It returns on success
a sequence of triples $(V_i,A_i',\mu_i)_{i=1}^N$, 
where $(V_i,\mu_i)$ is an optimal solution of the quantitative reach-avoid problem 
associated with $S$, $g$, $A_i'$ and the zero function. 
Moreover, $V_i(p) < \infty$ for all $p \in \cup_{i=1}^N A_j'$.
See \cite[Th.~V.1]{WeberKnoll21} for 
a proof for correctness of the algorithm in Fig.~\ref{a:coverage}. %
Key ingredient of this algorithm is the procedure 
\begin{equation}
\label{e:method:reachavoid}
\textsc{SolveReach}(S,g,A,H),
\end{equation}
which returns a pair of value function and optimal controller
for the reach-avoid problem associated with $S$, $g$, $A$ and $H$. 

\begin{figure}
\centering
\input{figures/symboliccontrol.tikz}
\caption{\label{fig:symboliccontrol}
Principle of symbolic controller synthesis \cite{ReissigWeberRungger17,WeberKnoll20}.}
\vspace{-3ex}
\end{figure}
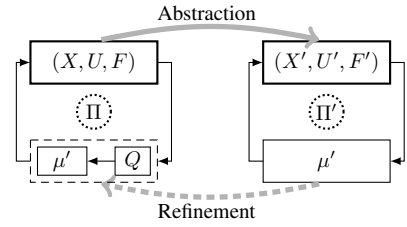

Computational techniques 
realizing this procedure in case of discrete data 
can be found in 
\cite{ReissigRungger18,
MacoveiciucReissig19,WeberKreuzerKnoll20,GirardEqtami21}.
In combination with symbolic controller synthesis \cite{ReissigWeberRungger17}
procedure \eqref{e:method:reachavoid} can be realized 
for sampled-data systems, 
whose underlying dynamics take the form \eqref{e:cont:dynamics}. 
Fig.~\ref{fig:symboliccontrol} depicts the principle. 
The given optimal control problem $\Pi$, 
which involves the sampled-data system $(X,U,F)$, 
is abstracted to an 
auxiliary optimal control problem $\Pi'$ 
that possesses \emph{finite} 
state and input space $X'$ and $U'$, respectively.
$X'$ is a finite cover of $X$ by non-empty sets 
and $U'$ a finite subset of $U$. 
$\Pi'$ is solved 
algorithmically. 
The resulting controller $\mu'$ is interconnected with 
a quantizer $Q \colon X \rightrightarrows X'$ 
and the composition $\mu'\circ Q$ solves $\Pi$ 
suboptimally on a predefined subset of $X$ \cite{ReissigRungger18}. 
Convergence results exist \cite{ReissigRungger18}. %
$Q$ is defined by $\Omega \in Q(x)$ iff $x\in \Omega$.
Applying this principle to coverage problems by means of Fig.~\ref{a:coverage}
results (on success) in controllers $\mu_1,\ldots,\mu_N$. 
The controller $\mu_i$ steers the state signal
near-optimally to $A_i$ whenever 
the state signal starts in $\cup_{i=1}^N A_i'$. 

This principle will be applied in our simulation results.
\section{Capacitated vehicle routing problem}
\label{s:vrp}
In this section, our contributions on the CVRP
are presented beginning with the problem definition, 
the proposed heuristic optimization, and 
concluding with simulation results.
Our contributions generalize 
those in \cite[Sect.~V]{WeberKnoll21}.
\subsection{Problem definition and heuristic optimization}
\label{ss:vrp:heuristic}
Firstly, we review the classical definition of a CVRP
and then we transfer it to our continuous domain.

In a classical \emph{capacitated} VRP, 
a depot is given,
$N - 1$ customer locations, distances between all customers and the depot, 
and a vehicle capacity $q$. 
A set of tours is sought such that 
1) all customers are served from the depot, 
2) a vehicle visits 
at most as many customers as its capacity, and
3) a given objective function is minimized along the tours. 

E.g., given three customers $2,3,4$ and depot $1$ with distance $1$ between each of them 
and vehicle capacity $2$ then 
the tours $1\to2\to1$ and $1\to 3 \to 4 \to 1$ constitute an optimal solution
minimizing the total distance travelled, which is $5$.

In what follows, we formalize the notions of tour, 
CVRP and its solutions in the 
variant that is considered in this work.
\begin{definition}
\label{def:tour}
Let $N \in \mathbb{Z}$, $N \geq 3$.
A finite sequence \[(1,t_2,\ldots,t_{N-1},1)\] where 
$t_i \in \mathbb{N}$, $t_i \geq 2$ for all $i \in \intco{2;N}$ and $t_i \neq t_j$ 
for all $i\neq j$ is called \term{$q$-tour} for each $q\geq N-2$.
\end{definition}
Therefore, a vehicle of capacity $q$ can perform at most a $q$-tour. 
The rigorous definition of the problem on a digraph with $N$ nodes and 
weighted adjacency matrix $C$ 
is as follows. 
\begin{definition}
\label{def:cvrp}
The triple
\begin{equation}
\label{e:cvrp}
(N,q,C) \in 
\mathbb{N} \times (\mathbb{N} \cup \{\infty\}) \times \mathbb{R}^{N\times N}
\end{equation}
is called a \term{classical (capacitated) vehicle routing problem} 
with vehicles of (uniform) capacity
$q$ and $N$ targets. 
A \term{feasible sequence of $q$-tours} to \eqref{e:cvrp} is sequence of $q$-tours
\begin{equation}
\label{e:cvrp:sol}
(\mathcal{T}_1,\ldots,\mathcal{T}_m)
\end{equation}
such that $\mathcal{T}_1,\ldots,\mathcal{T}_m \in \cup_{T=3}^{N+1}\intcc{1;N}^{\intcc{1;T}}$ and all elements of $\intoc{1;N}$
appear exactly once in \eqref{e:cvrp:sol}.
An \term{optimal solution} of \eqref{e:cvrp} is 
a feasible sequence of $q$-tours \eqref{e:cvrp:sol} that
minimizes 
\begin{equation}
\label{e:cvrp:objective}
\sum_{k=1}^m\nolimits \sum_{i=2}^{\operatorname{len} 
\mathcal{T}_k}\nolimits C_{\mathcal{T}_k(i-1),\mathcal{T}_k(i)}.
\end{equation}
\end{definition}

In the definition above, non-uniform customer demands are not considered and 
the number of tours/vehicles $m$ is part of the solution and not a constraint. 
However, the algorithm we are about to present immediately can be adopted
to other problem variants including, 
e.g., individual customer demands or 
the multiple travelling salesman problem \cite{BellmoreHong74}.
Also, other objectives than \eqref{e:cvrp:objective} can be considered. 

Definition \ref{def:cvrp} is transferred to 
the continuous domain as follows. 
Although it is
cumbersome to read, the idea is obvious. 
The constraint on the limited number of costumer supplies between two visits of the depot is modelled 
by means of $G_0$ in Def.~\ref{def:coverage} and using \eqref{e:cvrp:sol}.
\begin{definition}
\label{def:cvrp:continuous}
Let $q \in \mathbb{N}$. 
Let $\Pi$ be a coverage problem associated with $S$, $g$, $A_1,\ldots,A_N$
and $G_0$ defined by: 
$G_0(x|_{\intcc{0;t}}) = 0$ if there exist a feasible sequence of $q$-tours 
$(\mathcal{T}_1,\ldots,\mathcal{T}_m)$ 
such that for every 
$k,k' \in \intcc{1;m}$, $k \leq k'$ 
there exist $t_1 \leq \ldots \leq t_{\len \mathcal{T}_k} \leq t_2' \leq \ldots \leq t'_{\len \mathcal{T}_{k'}}$
such that 
\begin{equation*}
\forall_{l \in \intcc{1;\len \mathcal{T}_k}} : x(t_l) \in A_{\mathcal{T}_k(l)} \wedge 
\forall_{l \in \intoc{1;\len \mathcal{T}_{k'}}} : x(t'_l) \in A_{\mathcal{T}_{k'}(l)}.
\end{equation*}
Otherwise, $G_0(x|_{\intcc{0;t}}) = \infty$. 
Then $\Pi$ is called \term{capacitated vehicle routing problem} associated with
$S$, $g$, \term{target sets} $A_1,\ldots,A_N$ and \term{capacity} $q$. 
The set $A_1$ is called \term{depot}.
\end{definition}
Using the algorithmic methods of 
symbolic optimal control as depicted in Section \ref{s:review}, 
we can deduce a powerful heuristic to solve $\Pi$ 
in Def.~\ref{def:cvrp:continuous} suboptimally.
Algorithm \ref{a:cvrp} below 
provides controllers for each involved $q$-tour,
which are heuristically optimized
to reduce the cost functional \eqref{e:totalcost} as defined for $\Pi$. 
\begin{algorithm}
\caption{Capacitated vehicle routing problem}
\begin{algorithmic}[1]
\Input{$S$, $g$, $A_1, \ldots, A_N$, capacity $q$}
\State{\label{a:cvrp:solvecoverage}$(\bar V_i,A_i',\bar \mu_i)_{i=1}^N \gets \textsc{SolveCoverage}(S,g,A_1,\ldots,A_N)$}
\If{\label{a:cvrp:if}line \ref{a:cvrp:solvecoverage} failed}
\State{\Return{``Problem can't be solved"}}
\EndIf{\label{a:cvrp:endif}}
\State{\label{a:cvrp:matrix}Define $C \in \mathbb{R}^{N \times N}$ s. that $C_{i,j} = \min \{ \bar V_j(p) \mid p \in A_i'\}$}
\State{\label{a:cvrp:lkh}$(\mathcal{T}_1,\ldots,\mathcal{T}_m)\gets $ solution of classical VRP $(N,q,C)$}
\ForAll{\label{a:cvrp:for}$k \in \intcc{1;m}$}
\ForAll{$i \in \intco{1;\len\mathcal{T}_k}$}
\State{\label{a:cvrp:solvereach}\small{}$(V_i,\mu_{k,i}) \gets \textsc{SolveReach}(S,g,A_{\mathcal{T}_k(i)},\bar V_{\mathcal{T}_k(i+1)})$}
\EndFor{}
\EndFor{\label{a:cvrp:endfor}}
\Output{$\mu_{1,1},\ldots,\mu_{m,\operatorname{len} \mathcal{T}_m}$}
\end{algorithmic}
\label{a:cvrp}
\end{algorithm}
The algorithm works as follows.\\
\textit{\textbf{1st part (lines \ref{a:cvrp:solvecoverage}--\ref{a:cvrp:endif})}}. 
The coverage problem is solved. 
If successful (line \ref{a:cvrp:if}) 
the following steps will 
heuristically optimize the cost 
by computing suitable $q$-tours.\\
\textit{\textbf{2nd part (lines \ref{a:cvrp:matrix}--\ref{a:cvrp:lkh})}}. 
A classical CVRP is solved (line \ref{a:cvrp:lkh}), where
the value functions obtained before serve as cost estimates
between the sets $A_1,\ldots,A_N$ (line \ref{a:cvrp:matrix}). 
The involved minimum
is a heuristic choice and is relevant in the 3rd part.\\
\textit{\textbf{3rd part (lines \ref{a:cvrp:for}--\ref{a:cvrp:endfor})}}. 
Finally the actual
controller is computed. 
This last part builds upon the theoretical results in \cite[Sect.~III]{WeberKnoll20}
and is of significant importance to optimize the trajectory of each tour. 
In fact, in line \ref{a:cvrp:solvereach} an \emph{optimal} controller is computed 
for reaching \emph{firstly} $A_{\mathcal{T}_k(i)}$ and \emph{then} $A_{\mathcal{T}_k(i+1)}$.
Therefore, the minimum value in line \ref{a:cvrp:matrix} is a good estimate for the 
cost occurring between two targets.\\
\textit{\textbf{Output (and correctness)}}. 
The controller for the $q$-tour $\mathcal{T}_k$, $k \in \intcc{1;m}$
consists of the memoryless controllers $\mu_{k,1},\ldots,\mu_{k,\len \mathcal{T}_k}$ and 
is assembled as in Fig.~\ref{e:controller}. 
Specifically, the controller $\mu_{k,j}$ steers a
state signal starting in $\cup_{i=1}^N A'_i$ to $A_{\mathcal{T}_k(j)}$. 
This key property of Algorithm \ref{a:cvrp} is due to line \ref{a:cvrp:solvecoverage} 
since line \ref{a:cvrp:solvereach} does not change the suboptimality of $\mu_{k,j}$ \cite{WeberKnoll21}. 
\subsection{Experimental evaluation}
\label{ss:vrp:example}
\begin{figure}
\centering
\includegraphics[scale=.95]{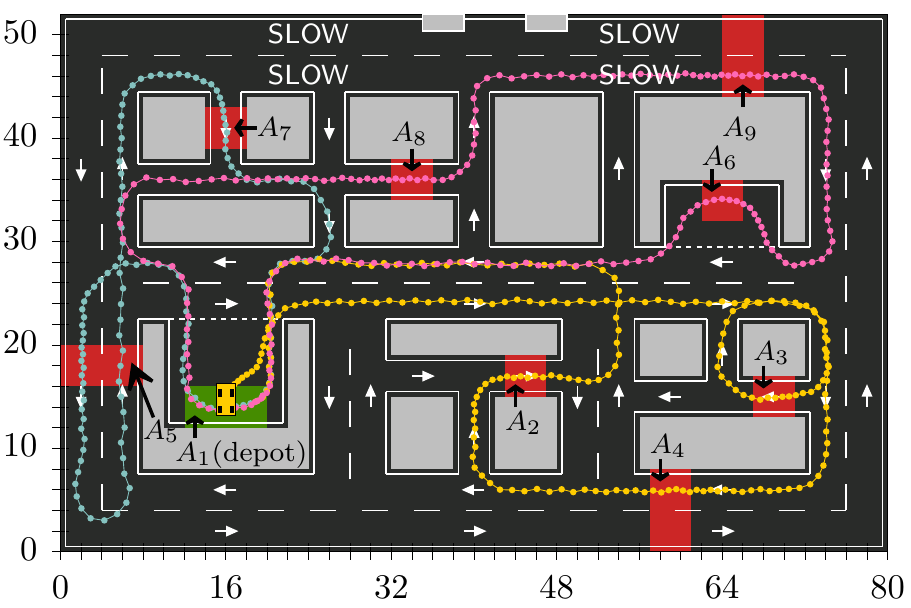}
\caption{\label{fig:cvrp}Truck delivery scenario of Section \ref{ss:vrp:example}.}
\vspace*{-\baselineskip}
\end{figure}

As an example for the application of Algorithm \ref{a:cvrp}, 
we consider a delivery truck 
in an urban environment as shown in Fig.~\ref{fig:cvrp}. 
(This scenario extends the one in \cite{WeberKnoll21}.)
The truck has capacity $3$ and shall visit 
the $8$ target areas $A_2,\ldots,A_9$ 
that are coloured in red.
The depot ($A_1$) is coloured in green. 
\subsubsection{Truck dynamics}
The truck is modelled 
in four dimensions according \eqref{e:cont:dynamics}
with $(x_1,x_2)$ being its position in the plane, 
and $x_3$ and $x_4$ its orientation and velocity, 
respectively. 
Control inputs
are the acceleration $u_1$ and the steering angle $u_2$. 
In fact, \eqref{e:cont:dynamics} is defined through $f\colon \R^4 \times U \to \R^4$,
\begin{align*}
&f(x,u) = \\ &(
x_4 \cos(\alpha + x_3) \cdot \beta,
x_4 \sin(\alpha + x_3) \cdot \beta, x_4 \tan(u_2),
u_1),
\end{align*}
where $U = \intcc{-6,4} \times \intcc{-0.5,0.5}$, $\alpha = \arctan(\tan(u_2)/2)$, 
$\beta =  \cos(\alpha)^{-1}$, $W=\{(0,0)\} \times \intcc{-0.01,0.01} \times \intcc{-0.1,0.1}$.
This control system is sampled with period $\tau = 0.1$
and a transition system $S$ in \eqref{e:system} is obtained 
with $X = \mathbb{R}^4$ and $U$ as above \cite[Def.~VIII.1]{ReissigWeberRungger17}. 
\subsubsection{Mission}
Formally, the capacitated vehicle routing problem $\Pi$
associated with $S$, $g$, target sets $A_1,\ldots, A_9$
and capacity $3$ is solved, 
where the involved objects are as follows. 
The coordinates of the target sets can be taken from Fig.~\ref{fig:cvrp}
and are similar to, e.g.,
$
A_6 = \intcc{62,66} \times \intcc{32,36} \times \mathbb{R} \times \intcc{0,7}.
$
The running cost $g$ satisfies $g(x,y,u) = \infty$ in three cases; 1)
if $x$ is outside the mission area
\begin{equation}
\label{e:cvrp:missionarea}
\intcc{0,80} \times \intcc{0, 52} \times \mathbb{R} \times \intcc{0,18},
\end{equation}
2) if $x$ is in one of the spacial obstacles depicted in Fig.~\ref{fig:cvrp} (gray coloured) or 3)
if $x$ violates the common right-hand traffic rules. 
E.g. the scenario contains a 
speed limit in the northernmost east-west road, so $x$ being in the set
\begin{equation*}
\intcc{0,80} \times \intcc{48,52} \times I_\mathrm{east} \times \intcc{0,12.5}, I_\mathrm{east} = \intcc{-\tfrac{3\pi}{8},\tfrac{3\pi}{8}}
\end{equation*}
makes $g$ finite.
If $g$ is finite, then $g$ balances minimum time and proper driving style. In fact,
\begin{equation}
g(x,y,u) = \tau + u_2^2 + \min_{m \in M} \| (y_1,y_2)  -  m \|_2,
\end{equation}
where $M \subseteq \intcc{0,30} \times \intcc{0,50}$. 
$M$ describes the axes of the roadways, 
e.g. $\intcc{2,50} \times \{2\} \subseteq M$, 
and proper approaches to the depot (see \cite[Sect.~VI.B]{WeberKnoll21}) and 
to target set $A_6$, e.g.
\begin{equation*}
\big ( \{60\} \times \intcc{28,34} \big ) \cup \big (\{68\} \times \intcc{28,34} \big ) \cup \big (\intcc{60,68} \times \{34\}\big ) \subseteq M.
\end{equation*}
\subsubsection{Application of Algorithm \ref{a:cvrp}}
Algorithm \ref{a:cvrp} is applied to an abstracted problem $\Pi'$ of $\Pi$ 
(cf. Fig.~\ref{fig:symboliccontrol}),
where lines \ref{a:cvrp:solvecoverage} and \ref{a:cvrp:solvereach}
are performed using value iteration \cite{WeberKreuzerKnoll20} and 
line \ref{a:cvrp:lkh} by \cite{Helsgaun17}. 
The problem $\Pi'$ contains
a discrete abstraction $(X',U',F')$ with 
$|X'| = 150 \cdot 114 \cdot 62 \cdot 50 \approx 53\cdot 10^6$, 
$|U'| = 9\cdot 11$. 
Here, the factors indicate how the abstraction is 
constructed on \eqref{e:cvrp:missionarea} and $U$ \cite{ReissigWeberRungger17}.
While the algorithm runs, 
$20$ reach-avoid problems are solved in $51$ minutes average. 
Total runtime is $17$ hours using $41$ GB RAM in parallel computation
with $26$ threads on Intel Xeon E5-2697 @ $2.60$ GHz.
\subsubsection{Obtained solution and evaluation}
The algorithm returns
the three tours 
$\mathcal{T}_1 = (1,3,4,2,1)$, 
$\mathcal{T}_2 = (1,5,7,1)$, 
$\mathcal{T}_3 = (1,8,9,6,1)$.
The cost \eqref{e:totalcost} for the mission beginning 
at the initial state $(16,14,\pi/2,0)$ is $244.6$. Fig.~\ref{fig:cvrp} 
shows the trajectory of this mission.
To rate the capabilities of our heuristic, 
we consider the tours $(1,2,4,1)$, $(1,5,7,8,1)$, $(1,9,3,6,1)$,
which are a ``reasonable" manual choice. 
Those tours have total cost $248.5$, 
which is $1.6\%$ more than the solution 
returned by Algorithm \ref{a:cvrp}. %
We would like to emphasize that 
in this example Algorithm \ref{a:cvrp} not only
returns a proper choice but also relieves the engineer of the task of picking tours among $8!=40320$
possible choices.
\section{Travelling salesman problem\\with arbitrary initial state} 
\label{s:tsp}
In this section, a technique for solving a \tsp{} is presented
where the initial state of the salesman is not necessarily 
inside the depot but can be somewhere in the state space. 
\subsection{Motivating example}
\label{ss:tsp:motivatingexample}
\begin{figure*}%
    \centering
    \subfloat[]{\includegraphics[scale=.54]{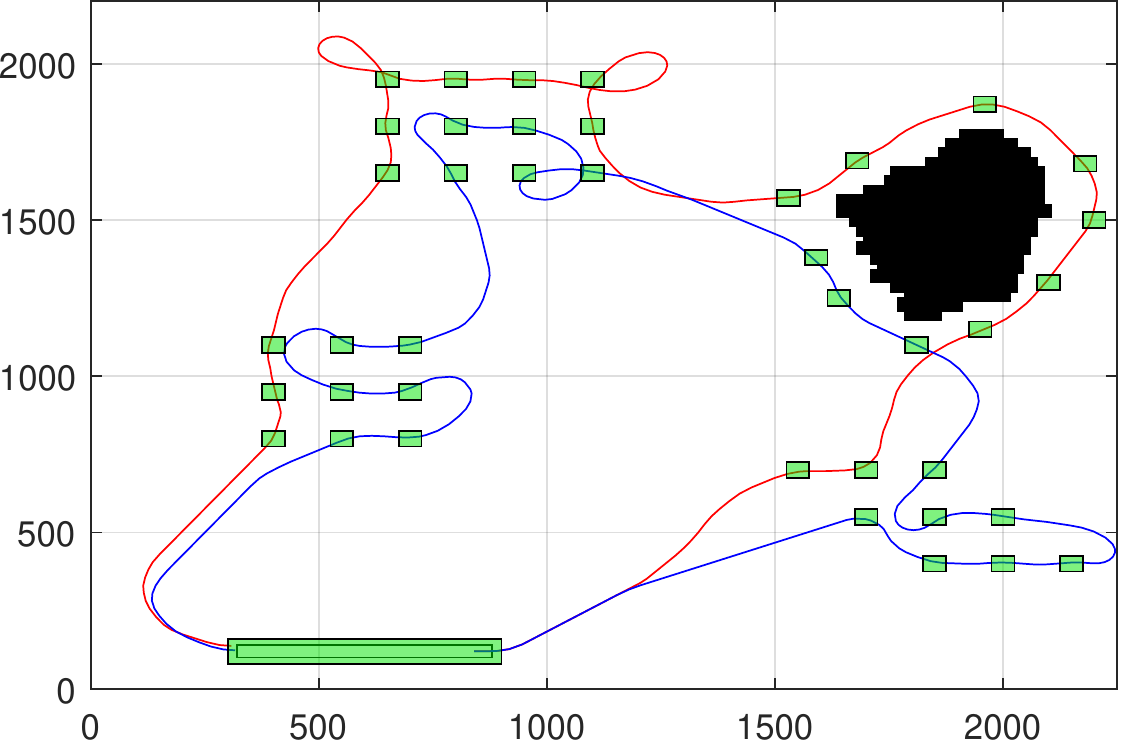}}%
    \quad
    \subfloat[]{\includegraphics[scale=.54]{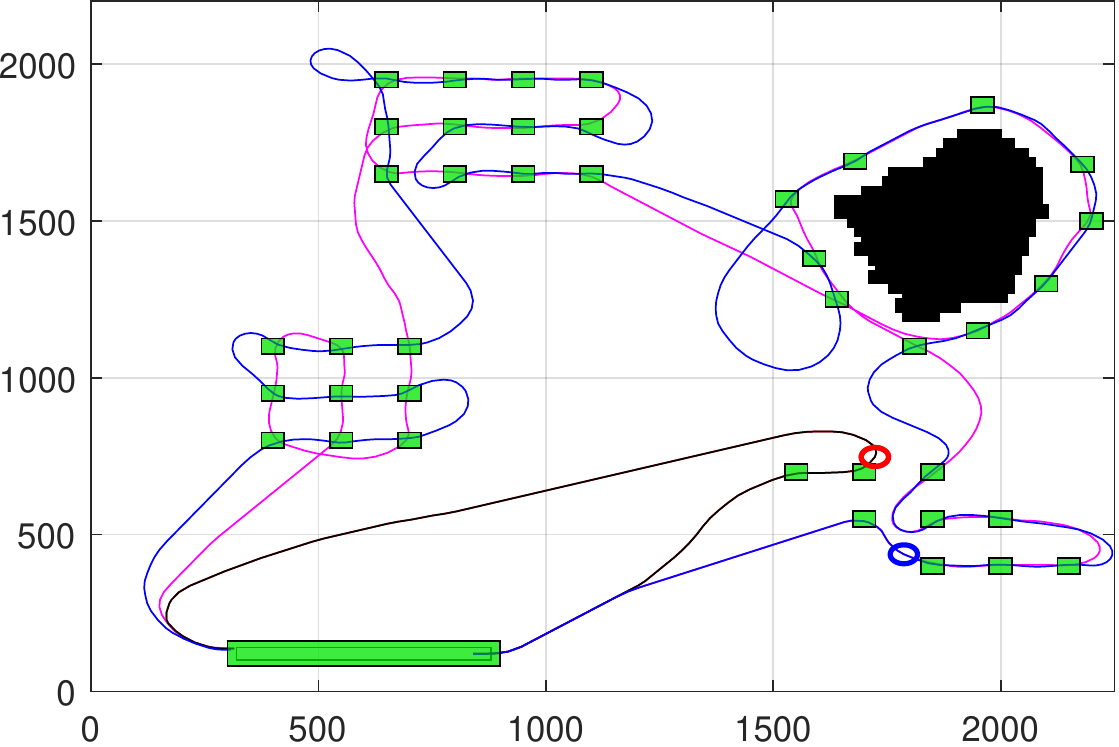}}%
    \quad
    \subfloat[]{%
    \begin{minipage}[b]{0.22\textwidth}
    \small         
    \begin{tabular}{l|r}
    Trajectory & Cost \\
    \hline \hline
    Red Fig.~\ref{fig:tsp}(a) & 193 \\
    Blue Fig.~\ref{fig:tsp}(a) & 222 \\
    \hline
    Blue Fig.~\ref{fig:tsp}(b) & 420 \\
    Magenta Fig.~\ref{fig:tsp}(b) & 382 \\
    \end{tabular}%
    \\[3ex]
    \small%
    \begin{tabular}{c|r}
    Lines Alg.~\ref{a:tsp} & Runtime \\
    \hline \hline
    \ref{a:tsp:solvecoverage} & 60 min. \\
    \ref{a:tsp:if} -- \ref{a:tsp:define} & $<$ 2 sec. \\
    \ref{a:tsp:line11} -- \ref{a:tsp:last} & 8 min. \\
    \ref{a:tsp:solvereach} with $\rho = 360$ & \textbf{10-15 sec.} 
    \end{tabular}%
    \end{minipage}%
}
    \caption{Reconnaissance mission with 2 UAVs in Section \ref{s:tsp}. The UAVs shall together visit 
    all of the 40 areas of interest (small green-coloured rectangles). The UAVs start from the airfield (green box near the bottom) and must return to it.
    (a) Trajectories of the two UAVs (coloured blue and red, respectively) in a mission \textbf{without} UAV failure.
    (b) Mission where the UAV coloured in red suffers damage at the red-circled state. 
    The UAV on the blue-coloured trajectory takes over the whole mission at the blue-circled state 
    according to the result of Algorithm \ref{a:tsp} (blue-coloured) and a comparison heuristic (magenta-coloured). 
    (c) Computational details to the missions depicted in (a) and (b). Implementation uses $24$ threads on Intel Xeon E5-2697 @ $2.60$ GHz.}%
    \label{fig:tsp}%
    \vspace*{-.5\baselineskip}
\end{figure*}

We would like to give a concrete 
example to which the algorithm 
to be presented can be applied. 
We consider a reconnaissance mission \cite{WeberKnoll21}, 
where two UAVs have to visit 40 areas of interest at minimum cost. 
The UAVs start and must return to the airfield, 
avoid obstacles, and 
share the task of visiting the areas. 
See Fig.~\ref{fig:tsp}(a).  
As stated, this is a VPR with a fixed number of vehicles 
of infinite capacity.
A controller for the mission can be obtained by Algorithm \ref{a:cvrp} 
(by adding the constraint of a fixed number of vehicles to line \ref{a:cvrp:lkh}).
Now, we assume that at some random point in time during mission 
one UAV suffers damage and must return to base right away.
The remaining UAV shall fly to all 
not yet visited areas of interest 
including those initially assigned to the other UAV. 
See Fig.~\ref{fig:tsp}(b).
In this case, a TSP must be solved with initial state 
the current state of the functioning UAV.
And clearly, the problem must be solved in short time 
so that the UAV can continue its mission without notable delay.\looseness=-1
\subsection{Problem definition and heuristic optimization}
\label{s:tsp:heuristics}
Formally, we consider the following problem.
\begin{definition}
Let $\Pi$ be a coverage problem associated with 
$A_1,\ldots,A_N \subseteq X$ and $G_0$
defined by $G_0(x_{\intcc{0;t}}) = \infty$ 
if $x(t) \notin A_1$ and otherwise $G_0(x_{\intcc{0;t}})=0$. 
Then $\Pi$ is called \term{travelling salesman problem} with 
\term{target sets $A_1,\ldots,A_N$}. 
The set $A_1$ is called \term{depot}.
\end{definition}
The algorithm below provides a suboptimal solution of 
the \tsp{} that is heuristically optimized for trajectories 
starting at the fixed initial state $x_0$.
Moreover, it includes a technique to speed up computing time. 
The structure is similar to the one of Algorithm \ref{a:cvrp}:\\
\begin{algorithm}
\caption{TSP with arbitrary initial state}
\begin{algorithmic}[1]
\Input{$S = (X,U,F)$, $g$, $A_1, \ldots, A_N$, initial state $x_0 \in X$}
\Parameter{$\rho \in \mathbb{R}_+$}
\State{\label{a:tsp:solvecoverage}$(V_i,A_i',\mu_i)_{i=1}^N \gets \textsc{SolveCoverage}(S,g,A_1,\ldots,A_N)$}
\If{\label{a:tsp:if}line \ref{a:cvrp:solvecoverage} failed \textbf{or} $V_1(x_0) = \infty$}
\State{\Return{``Problem can't be solved"}}
\EndIf{\label{a:tsp:endif}}
\ForAll{\label{a:tsp:line2}$i,j \in \{1,\ldots,N\}$}
\If{$i = 1 \text{ and } j > 1$}
\State{$C_{1,j} \gets V_j(x_0)$}
\Else{}
\State{$C_{i,j} \gets \min\{V_j(p) \mid p \in A_i'\}$}
\EndIf{}
\EndFor{\label{a:tsp:line8}}
\State{\label{a:tsp:lkh}$\mathcal{T} \gets $ solution of classical VRP $(N,\infty,C)|_{\text{no. of vehicles}=1}$}
\State{\label{a:tsp:define}Define $h_i \colon X \times X \times U\to \mathbb{R}_+ \cup \{\infty\}$, $i \in \intcc{1;N}$ by:
\begin{equation*}
h_i(x,y,u) = \begin{cases}
g(x,y,u) , & \text{if } \forall_{j \in \intcc{1;n}} \inf_{a \in A_i}\limits | x_j - a_j| \leq \rho\\ 
\infty , & \text{otherwise}
\end{cases}
\end{equation*}
}
\ForAll{\label{a:tsp:line11}$i \in \intco{1;N}$}
\State{\label{a:tsp:solvereach}$(\bar V_i,\bar \mu_i) \gets \textsc{SolveReach}(S,h_i,A_{\mathcal{T}(i)},V_{\mathcal{T}(i+1)})$}
\ForAll{$x \in X$}
\If{$\bar V_i(x) = \infty$}
\State{$\bar \mu_i(x) = \mu_i(x)$}
\EndIf{}
\EndFor{}
\EndFor{\label{a:tsp:last}}
\Output{$\bar \mu_1, \ldots, \bar \mu_N$}
\end{algorithmic}
\label{a:tsp}
\end{algorithm}%
\textit{\textbf{1st part (lines \ref{a:tsp:solvecoverage}--\ref{a:tsp:endif})}}.
The coverage problem is solved (line \ref{a:tsp:solvecoverage}). 
If successful (including $x_0$), the algorithm continues.\\
\textit{\textbf{2nd part (lines \ref{a:tsp:line2}--\ref{a:tsp:lkh})}}. 
A classical \tsp{} is solved, where the underlying weighted adjacency matrix $C$
is roughly defined as follows.
In the first row the estimated costs 
for the initial state $x_0$ to the targets $A_2,\dots,A_N$ appear 
(lines \ref{a:tsp:line2}--\ref{a:tsp:line8}), 
where the previously computed value functions serve as cost estimates.
Otherwise, the best-possible cost from $A_i$ to $A_j$ is contained in $C_{i,j}$ \cite{WeberKnoll20}.\\
\textit{\textbf{3rd part (lines \ref{a:tsp:define}--\ref{a:tsp:last})}}. 
The actual controllers are computed.
As seen in Section \ref{ss:vrp:example}, 
this is time consuming if all the occurring
reach-avoid problems (line \ref{a:tsp:solvereach}) 
are fully solved. 
Therefore, in Algorithm \ref{a:tsp} the reach-avoid problems 
are only solved on the neighbourhood of $A_i$, 
which is quantified by the parameter $\rho$. 
Finally, the controllers are combined with the already computed 
controllers of the coverage specification.
In this way, line \ref{a:tsp:solvereach} can be executed 
while the closed loop is active at the sacrifice of further optimality.\\
\textit{\textbf{Output (and correctness)}}. Controller $\bar \mu_i$ steers a state signal stating in $\{x_0\} \cup A_1 \cup \ldots \cup A_N$ to $A_i$ by line \ref{a:tsp:solvecoverage}
since the rest of the algorithm does not change the suboptimality of $\bar \mu_i$.
\subsection{Experimental evaluation}
Two UAV scenarios were simulated 
to demonstrate the effectiveness of our heuristic. 
The first scenario has already been described 
in Section \ref{ss:tsp:motivatingexample}. 
The second scenario, depicted in Fig.~\ref{fig:prunus}, 
has larger expansion and involves ten UAVs. 
\subsubsection{UAV dynamics}
The UAVs are modelled by Dubins dynamics \cite{LaValle06}
given by \eqref{e:cont:dynamics} with $f\colon \mathbb{R}^3 \times U \to \mathbb{R}^3$, 
$U = \intcc{20,50} \times \intcc{-0.5,0.5}$,
$f(x,u) = (u_1 \cos(x_3) , u_1 \sin(x_3), u_2)$,
where the pair $(x_1,x_2)$ locates the vehicle in the plane and $x_3$ is the 
course angle. 
Control inputs are the velocity $u_1$ and the course change rate $u_2$. 
Wind in the mission area 
is assumed and taken into account by 
$W = \intcc{-5,5} \times \intcc{-2,2} \times \intcc{-0.04,0.04}$ in \eqref{e:cont:dynamics}.
The UAV dynamics is sampled with period $\tau =  0.65$ to obtain
a sampled system $S$ with dynamics \eqref{e:discrete:dynamics} \cite[Def.~VIII.1]{ReissigWeberRungger17}.
\subsubsection{Running cost}
\begin{figure*}
\subfloat[]{
\includegraphics[scale=.7]{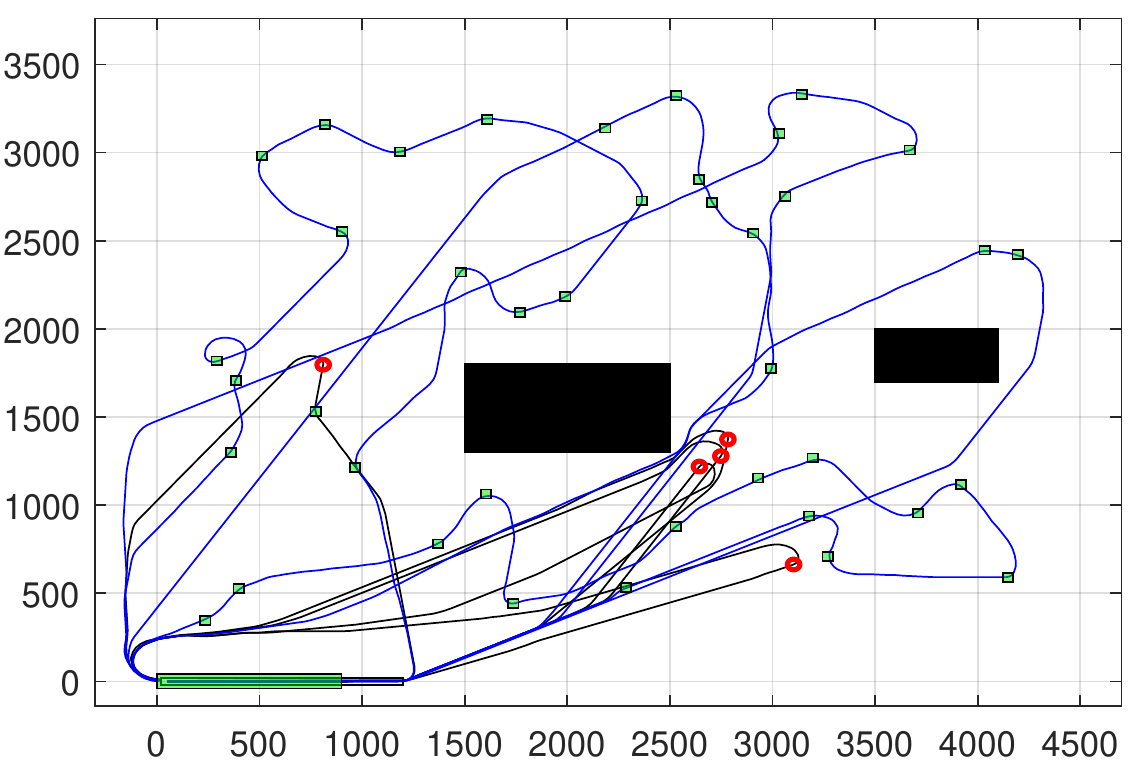}}\qquad 
\subfloat[]{\includegraphics[scale=.7]{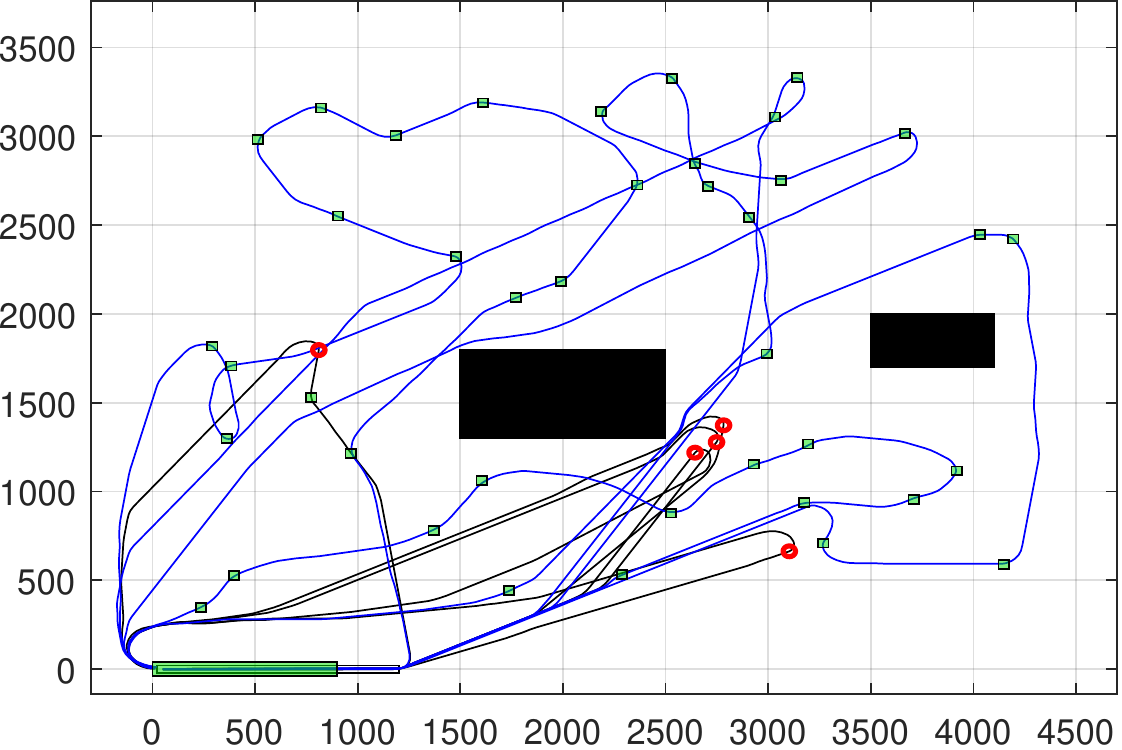}}
\caption{\label{fig:prunus}UAV mission with 10 UAVs and randomly distributed areas of interest (small green-coloured rectangles) in Section \ref{s:tsp}. (a) Mission guidance by Algorithm \ref{a:tsp}. 
(b) Mission guidance by the comparison heuristic.
}
\vspace*{-\baselineskip}
\end{figure*}

The running cost for both missions balances 
minimum time and small course rate change, 
and includes the requirement to avoid obstacles, 
i.e. $g$ in \eqref{e:runningcost} is given by
\begin{equation*}
g(x,y,u) = \begin{cases}
\infty, & \text{if } x \in (\R^3 \setminus X_{\mathrm{mis}}) \cup A_{\mathrm{nofly}} \cup A_{\mathrm{hill}} \\
\tau + u_2^2, & \text{otherwise (} u_2 \text{ in radians)}
\end{cases}
\end{equation*}
where
$X_\mathrm{mis} = \intcc{0,2500} \times \intcc{0,2200} \times \mathbb{R}$
for the mission with two UAVs and 
$X_\mathrm{mis} = \intcc{-300,4700} \times \intcc{-140,3760} \times \mathbb{R}$
for the mission with ten UAVs. 
The spatial obstacle set $A_{\mathrm{hill}}$ is given by the 
black-coloured regions in Fig.~\ref{fig:tsp} and Fig.~\ref{fig:prunus}, respectively. 
The set $A_{\mathrm{nofly}}$ given by $\intcc{320,880} \times \intcc{100,140} \times \intcc{12^\circ,348^\circ}$ and 
by $\intcc{20,1200} \times \intcc{-20,20} \times \intcc{12^\circ,348^\circ}$, 
respectively, models a proper air traffic guidances for takeoff and landing. 
\subsubsection{Comparison heuristics}
\label{sss:comparison}
To validate the efficiency of our heuristic, 
the results of Algorithm \ref{a:tsp} are  
compared with another heuristic method as follows.
Let $\mu_1,\ldots,\mu_{41}$ and 
$V_1,\ldots,V_{41}$ be as in line \ref{a:tsp:solvecoverage} of 
Algorithm~\ref{a:tsp}.
At the moment of failure of one UAV, 
its target areas are transferred to 
the UAV that has the shortest Euclidean distance 
to the failure location at the time of the failure.
This vehicle flies to its assigned target areas 
in the further course depending on the current cost, i.e. 
selects at the moment of the stopping of the current controller at state $x$  
the controller $\mu_i$ according to
\begin{equation*}
i \in \argmin \{ V_i(x) \mid i \in \intcc{2;41}\}.
\end{equation*}
Thus, after a target is reached, 
values of the value functions to all other remaining target areas 
(airfield is excluded) are compared. 
The target with the lowest value of the value function 
is selected as next target.
\subsubsection{Scenario with two UAVs}
We consider the case where the UAV on the red-coloured trajectory in Fig.~\ref{fig:tsp}(a) 
must abruptly return to the airfield when being in the red-circled state in Fig.~\ref{fig:tsp}(b). 
At that time, the UAV on the blue-coloured trajectory is 
in the blue-circled state in Fig.~\ref{fig:tsp}(b)
and shall take over all tasks.
To this end, Algorithm \ref{a:tsp} is applied to $S$, $g$, the remaining $38$ target sets and 
the indicated initial state. 
The parameter $\rho$ is chosen as $360$. 
The resulting trajectory is depicted in Fig.~\ref{fig:tsp}. 
If lines \ref{a:tsp:define}--\ref{a:tsp:last} in Algorithm \ref{a:tsp} are skipped
the resulting cost of the trajectory is $8\%$ higher 
than the optimized one.
Nevertheless, in this scenario 
our heuristic proved weaker than the method based on 
the frequent evaluation of the value function. 
Runtimes and costs are given in Fig.~\ref{fig:tsp}(c).
\subsubsection{Scenario with 10 UAVs}
While the number of target areas remains the same in this scenario, 
the number of UAVs increases to ten and 
the target sets are arranged much more diffuse. See Fig.~\ref{fig:prunus}. 
The error-free mission is computed by means of a VRP with 10 vehicles as constraint and infinite capacity.
We let five UAVs be damaged at the red-circled states in Fig.~\ref{fig:prunus}
and apply Algorithm \ref{a:tsp} and 
compare the results with the heuristic in Subsection \ref{sss:comparison}.
In this scenario Algorithm \ref{a:tsp} proves more efficient. It provides
controllers leading to a mission cost of $1994$ while for the naive solution it is $2023$. 
\section{Conclusions}
\label{s:conclusions}
The application of symbolic optimal control
to vehicle mission guidance has been 
demonstrated for two types of missions. 
Following \cite{WeberKnoll21} 
we showed how to utilize the value functions
of the underlying reach-avoid tasks and classical solvers for VRPs 
in order to heuristically reduce the overall mission cost (Section \ref{s:vrp}).
Although symbolic optimal control is currently not suitable
for missions with spontaneous changes in the specification due to high runtime requirements,
we have given a possible 
direction towards dynamic controller adaptation (Section \ref{s:tsp}).
\bibliography{../../BibTex/mrabbrev,../../BibTex/articles,../../BibTex/books,../../BibTex/online}
\end{document}

%% file: figures/highlevelloop.tikz
\newcommand{\xca}{-3.}
\newcommand{\yca}{-.3}
\newcommand{\ycb}{.7}
\newcommand{\xcb}{3.}
\newcommand{\xsa}{\xca}
\newcommand{\ysa}{1.5}
\newcommand{\xsb}{\xcb}
\newcommand{\ysb}{2.9}
\newcommand{\xma}{-8.2}
\newcommand{\xmb}{-7.0}
\newcommand{\ymb}{7}
\newcommand{\ofs}{1.4} 
\newcommand{\ofss}{-5.1} 
\newcommand{\ofsw}{1.4} 
\begin{tikzpicture}[scale=0.55,>=latex]
\draw[thick] (\xma,\ysa+\ofss) rectangle node[rotate=90,text width=3cm,align=center]{scheduler} (\xmb,\ymb+\ofss-.2);
\draw[thick] (\xsa,\ysa) rectangle node[text width=3cm,align=center]{Plant $S$} (\xsb,\ysb);
\draw[thick] (\xca,\yca) rectangle node[text width=3cm,align=center]{Controller $\mu_1$} (\xcb,\ycb);
\draw[thick,fill=black!15!white] (\xca,\yca-\ofs) rectangle node[text width=3cm,align=center]{Controller $\mu_2$} (\xcb,\ycb-\ofs);
\draw[thick] (\xca,\yca-2.7*\ofs) rectangle node[text width=3cm,align=center]{Controller $\mu_K$} (\xcb,\ycb-2.7*\ofs);
\coordinate (v1) at (\xsb,\ysa/2+\ysb/2) ;
\coordinate (v2) at (\xsb+\ofsw,\ysa/2+\ysb/2) ;
\coordinate (v3) at (\xsb+\ofsw,\yca/2+\ycb/2) ;
\coordinate (v4) at (\xsb,\yca/2+\ycb/2) ;
\coordinate (v5) at (\xsb+\ofsw,\yca/2+\ycb/2 - \ofs) ;
\coordinate (v6) at (\xsb,\yca/2+\ycb/2 - \ofs) ;
\coordinate (v8) at (\xsb+\ofsw,\yca/2+\ycb/2 - 2.*\ofs) ;
\coordinate (v9) at (\xsb+\ofsw,\yca/2+\ycb/2 - 2.7*\ofs) ;
\coordinate (v10) at (\xsb,\yca/2+\ycb/2 - 2.7*\ofs) ;
\coordinate (v11) at (\xsb,\ysa/2+\ysb/2+\ofss) ;
\coordinate (v12) at (\xsa,\yca/2+\ycb/2);
\coordinate [draw,circle,scale=0.4] (v13) at (\xsa-\ofsw,\yca/2+\ycb/2);
\coordinate (v14) at (\xsa,\yca/2+\ycb/2-\ofs);
\node [draw,circle,scale=0.4] (v15) at (\xsa-\ofsw,\yca/2+\ycb/2-\ofs) {};
\coordinate (v16) at (\xsa,\yca/2+\ycb/2-2.7*\ofs) ;
\coordinate [draw,circle,scale=0.4] (v17) at (\xsa-\ofsw,\yca/2+\ycb/2-2.7*\ofs);
\node [draw,circle,scale=0.4] (v18) at (-5.75,-2) {};
\coordinate (v20) at (\xsa,\ysa/2+\ysb/2);
\coordinate (v21) at (-5,-1.75) {} {}; 
\coordinate (v22) at (\xma/2+\xmb/2,-3.9) {} ; 
\coordinate (v23) at (-6.,\ysa/2+\ysb/2);
\coordinate (v24) at (\xma/2+\xmb/2,\ymb+\ofss-.2);
\coordinate (v25) at (-5.75,\ysa/2+\ysb/2);
\coordinate (v26) at (\xma/2+\xmb/2,\ysa+\ofss);
\coordinate (v27) at (-6.25,\ysa/2+\ysb/2) {};
\coordinate (v28) at (-6.25,-2.05) {};
\coordinate (v29) at (-6.1,\ysa/2+\ysb/2) {};
\coordinate (v30) at (-6.1,-1.95) {};
\coordinate (v31) at (-5.85,-2.05);
\coordinate (v32) at (-5.85,-1.95);
\coordinate (v33) at (\xsa-\ofsw+.2,\yca/2+\ycb/2);
\coordinate (v34) at (\xsa-\ofsw+.1,\yca/2+\ycb/2-2.7*\ofs+.1);
\coordinate (v35) at (\xsa-\ofsw,\yca/2+\ycb/2-\ofs-.15);
\draw[<-,thick] (v24) |- (v27) -- (v28) -- (v31);
\draw[->] (v32) -- (v30) |- (v20);
\draw[->] (v1) -- (v2) -- (v3) -- (v4);
\draw[->] (v3) -- (v5) -- (v6);
\draw[->] (v8) -- (v9) -- (v10);
\draw[-] (v12) -- (v13);
\draw[-] (v14) -- (v15);
\draw[-] (v16) -- (v17);
\draw[thick,-] (v18) -- (v35);
\draw[dashed,thick,-] (v18) -- (v33);
\draw[dashed,thick,-] (v18) -- (v34);
\draw[thick,->] (v26) -- (v22) -| (v21); 
\node at (\xsa/2+\xsb/2,\yca/2+\ycb/2-1.75*\ofs) {$\vdots$};
\node at (-3.75,\yca/2+\ycb/2-1.75*\ofs) {$\vdots$};
\node at (\xsb+\ofsw,-1.75) {$\vdots$};
\node at (3.5,\ysb-.4) {$x$};
\node at (-4.5,\ysb-.4) {$u$};
\node at (-7,\ysb-.4) {$v$};
\node at (-3.75,0.6) {\footnotesize$(u,v)$};
\node at (-3.75,0.6-\ofs) {\footnotesize$(u,v)$};
\node at (-3.75,0.6-5*\ofs/2-.35) {\footnotesize$(u,v)$};
\end{tikzpicture}

%% file: figures/symboliccontrol.tikz
\usetikzlibrary{arrows}
\begin{tikzpicture}[scale=.75, every node/.style={scale=.8}]
\draw[thick]  (-3.75,3.25) rectangle (-1.5,2.5) node[pos=.5] {$(X,U,F)$};
\draw[thick]   (0.375,3.25) node (v2) {} rectangle (2.625,2.5) node[pos=.5] {$(X',U',F')$};
\draw[-latex] (-1.5,2.875) -- (-1.25,2.875) -- (-1.25,1.125) -- (-1.5,1.125);
\draw  (-2.25,1.375) rectangle (-1.625,0.875) node [pos=.5] (v4) {$Q$};
\draw  (-3.625,1.375) rectangle (-2.75,0.875) node[pos=.5] {$\mu'$};
\draw[thick,densely dotted] (-2.625,2) circle (0.3) node {$\Pi$};
\draw[thick,densely dotted] (1.5,2) circle (0.3) node {$\Pi'$};
\draw[-latex] (-2.25,1.125) -- (-2.75,1.125);
\draw[-latex] (-3.75,1.125) -- (-4,1.125) -- (-4,2.875) -- (-3.75,2.875);
\node (v1) at (-2.625,3.25) {};
\node (v4) at (1.5,3.25) {};
\draw[gray!60!white,line width=2][->]  (v1) edge [ bend angle=12,bend left] (v4);
\node at (-0.625,3.75) {Abstraction};
\draw  (0.375,1.5) rectangle (2.625,0.75) node[pos=.5] {$\mu'$};
\draw[-latex] (2.625,2.875) -- (2.875,2.875) -- (2.875,1.125) -- (2.625,1.125);
\draw[-latex] (0.375,1.125) -- (0.125,1.125) -- (0.125,2.875) -- (0.375,2.875);
\node (v3) at (0.375,0.75) {};
\node (v5) at (1.5,0.75) {};
\node (v6) at (-2.625,0.75) {};
\draw[gray!60!white,line width=2,densely dashed][->]  (v5) edge [ bend angle=12,bend left] (v6);
\node at (-0.625,0.25) {Refinement};
\node (v7) at (-1.5,0.75) {};
\draw[densely dashed]  (-3.75,1.5) rectangle (v7);
\end{tikzpicture}